\documentclass[11pt]{article}
\usepackage[dvips]{graphicx}
\usepackage{color}
\usepackage{amsmath}
\usepackage{amsfonts}
\usepackage{amssymb}
\usepackage{amsthm}
\usepackage{newlfont}
\usepackage{epsfig}
\usepackage{multirow}
\usepackage{setspace}
\usepackage{hyperref}
\usepackage{cite}

\begin{document}

\newtheorem{assumption}{Assumption}[section]
\newtheorem{definition}{Definition}[section]
\newtheorem{lemma}{Lemma}[section]
\newtheorem{proposition}{Proposition}[section]
\newtheorem{theorem}{Theorem}[section]
\newtheorem{corollary}{Corollary}[section]
\newtheorem{remark}{Remark}[section]

\title{A Linear Programming Approach to Weak Reversibility and Linear Conjugacy of Chemical Reaction Networks}
\author{Matthew D. Johnston$^a$, David Siegel$^a$ and G\'{a}bor Szederk\'{e}nyi$^{b,c}$ \bigskip \\
${}^a$ Department of Applied Mathematics,\\
University of Waterloo,\\
Waterloo, Ontario, Canada N2L 3G1\\
${}^b$ (Bio)Process Engineering Group, IIM-CSIC,\\ Spanish National Research Council,\\ C/Eduardo Cabello, 6, 36208 Vigo, Spain \\
${}^c$ Computer and Automation Research Institute,\\
Hungarian Academy of Sciences\\
H-1518, P.O. Box 63, Budapest, Hungary}
\date{}
\maketitle

\tableofcontents

\bigskip

\begin{abstract}

A numerically effective procedure for determining weakly reversible chemical reaction networks that are linearly conjugate to a known reaction network is proposed in this paper. The method is based on translating the structural and algebraic characteristics of weak reversibility to logical statements and solving the obtained set of linear (in)equalities in the framework of mixed integer linear programming. The unknowns in the problem are the reaction rate coefficients and the parameters of the linear conjugacy transformation. The efficacy of the approach is shown through numerical examples.

\end{abstract}


\bigskip

\section{Introduction}

A chemical reaction network is given by sets of reactants reacting at prescribed rates to form sets of products. The mathematical study of such networks has been applied recently to such fields as industrial chemistry, systems biology, gene regulation, among others \cite{C-P-F,H2,S}. There has also been significant theoretical work in the literature on such questions as persistence \cite{A,A-S,A3}, multistability \cite{C-F1, C-F2,S-F}, monotonicity \cite{B-A,B2}, the global attractor conjecture for complex balanced systems \cite{A,A-S,C-D-S-S,C-N-P}, and conjugacy of reaction networks \cite{C-P,J-S2}.

One problem which has attracted recent attention has been that of determining when two reaction networks can exhibit the same qualitative dynamics despite disparate network structure. It has been long known that two networks can given rise to the same governing set of differential equations under the assumption of mass-action kinetics \cite{K,H-J1}. In \cite{C-P} and \cite{Sz1}, the authors complete the question of what network structures can given rise to a set of governing differential equations. This work is extended in \cite{J-S2} to networks which do not necessarily have the same set of differential equations but rather have trajectories related by a non-trivial linear transformation.

The problem of algorithmically determining when a network has the same governing dynamics as another network satisfying specified properties was first addressed in \cite{Sz2} where the author presented a mixed-integer linear programming (MILP) algorithm capable of determining sparse and dense realizations, i.e. networks with the fewest or greatest number of reactions capable of generating the given dynamics. This line of research was continued in \cite{Sz-H} in which the authors extended the algorithm to complex and detailed balanced networks, and in \cite{Sz-H-T} in which the authors addressed the problem of finding dense weakly reversible realizations.

In this paper we show how the problem of determining weakly reversible realizations presented in \cite{Sz-H-T} can be reformulated as a linear constraint within the established MILP framework. This reformulation significantly eases the computational cost associated with the problem which had previously been solved through successive applications of an optimization algorithm and checking for weak reversibility with one of Kosaraju's, Tarjan's or Gabow's algorithm \cite{BJ-G,N-SS}. We also extend the algorithm to encompass the linearly conjugate networks introduced in \cite{J-S2}. We show how the algorithm can be used to effortlessly reproduce results from the literature and easily handle large-scale networks not yet considered.

\section{Background}

In this section we present the terminology and notation relevant to chemical reaction networks and the main results from the literature upon which we will be building.

\subsection{Chemical Reaction Networks}

We will consider the chemical \emph{species} or \emph{reactants} of a network to be given by the set $\mathcal{S} = \left\{ X_1, X_2, \ldots, X_n \right\}$. The combined elements on the left-hand and right-hand side of a reaction are given by linear combinations of these species. These combined terms are called \emph{complexes} and will be denoted by the set $\mathcal{C} = \left\{ C_1, C_2, \ldots, C_m \right\}$ where
\[C_i = \sum_{j=1}^n \alpha_{ij} X_j, \; \; \; i=1, \ldots, m\]
and the $\alpha_{ij}$ are nonnegative integers called the \emph{stoichiometric coefficients}. We define the reaction set to be $\mathcal{R} = \left\{ (C_i,C_j) \; | \; C_i \mbox{ reacts to form } C_j \right\}$ where the property $(C_i,C_j) \in \mathcal{R}$ will more commonly be denoted $C_i \to C_j$. To each $(C_i,C_j) \in \mathcal{R}$ we will associate a positive \emph{rate constant} $k_{ij} > 0$ and to each $(C_i,C_j) \not\in \mathcal{R}$ we will set $k_{ij} = 0$. The triplet $\mathcal{N} = (\mathcal{S}, \mathcal{C}, \mathcal{R})$ will be called the \emph{chemical reaction network}.

The above formulation naturally gives rise to a directed graph $G(V,E)$ where the set of vertices is given by $V = \mathcal{C}$, the set of directed edges is given by $E = \mathcal{R}$, and the rate constants $k_{ij}$ corresponds to the weights of the edges from $C_i$ to $C_j$. In the literature this has been termed the \emph{reaction graph} of the network \cite{H-J1}. Since complexes may be involved in more than one reaction, as a product or a reactant, there is further graph theory we may consider. A \emph{linkage class} is a maximally connected set of complexes, that is to say, two complexes are in the same linkage class if and only if there is a sequence of reactions in the reaction graph (of either direction) which connects them. A reaction network is called \emph{reversible} if $C_i \to C_j$ for any $C_i, C_j \in \mathcal{C}$ implies $C_j \to C_i$. A reaction network is called \emph{weakly reversible} if $C_i \to C_j$ for any $C_i, C_j \in \mathcal{C}$ implies there is some sequence of complexes such that $C_i = C_{\mu(1)} \to C_{\mu(2)} \to \cdots \to C_{\mu(l-1)} \to C_{\mu(l)} = C_j$.

A directed graph is called \emph{strongly connected} if there exists a directed path from each vertex to every other vertex. A \emph{strongly connected component} of a directed graph is any set of vertices for which paths exists from each vertex in the set to every other vertex in the set. For a weakly reversible network, the linkage classes clearly correspond to the strongly connected components of the reaction graph.

Assuming mass-action kinetics, the dynamics of the specie concentrations over time is governed by the set of differential equations
\begin{equation}
\label{de}
\frac{d\mathbf{x}}{dt} = Y \cdot A_k \cdot \Psi(\mathbf{x})
\end{equation}
where $\mathbf{x} = [ x_1 \; x_2 \; \cdots \; x_n ]^T$ is the vector of reactant concentrations. The \emph{stoichiometric matrix} $Y$ contains entries $[Y]_{ij} = \alpha_{ji}$ and the \emph{Kirchhoff} or \emph{kinetics} matrix $A_k$ is given by
\begin{equation}
\label{kinetics}
[A_k]_{ij} = \left\{ \begin{array}{cll} -\sum_{l=1,l \not= i}^m k_{il}, & \mbox{  if  } & i = j \\ k_{ji} & \mbox{  if  } & i \not= j. \end{array} \right.
\end{equation}
Finally, the \emph{mass-action vector} $\Psi(\mathbf{x})$ is given by
\begin{equation}
\label{psi}
\Psi_j(\mathbf{x}) = \prod_{i=1}^n x_i^{[Y]_{ij}}, \; \; \; j=1, \ldots, m.
\end{equation}

\subsection{Sparse and Dense Realizations}

Under the assumption of mass-action kinetics, it is possible for two reaction networks $\mathcal{N}^{(1)}$ and $\mathcal{N}^{(2)}$ to give rise to the same set of governing differential equations. In other words, it is possible that
\[Y^{(1)} \cdot A_k^{(1)} \cdot \Psi^{(1)}(\mathbf{x}) = Y^{(2)} \cdot A_k^{(2)} \cdot \Psi^{(2)}(\mathbf{x}) = \mathbf{f}(\mathbf{x}), \; \; \; \forall \; \mathbf{x} \in \mathbb{R}_{>0}^n\]
where $Y^{(i)}$, $A_k^{(i)}$, and $\Psi^{(i)}(\mathbf{x})$, $i=1,2$, are the stoichiometric matrices, kinetics matrices, and mass-action vectors defined for $\mathcal{N}^{(1)}$ and $\mathcal{N}^{(2)}$ respectively. The networks $\mathcal{N}^{(1)}$ and $\mathcal{N}^{(2)}$ are called different \emph{realizations} of the kinetics $\mathbf{f}(\mathbf{x})$ although it will sometimes be more convenient to consider $\mathcal{N}^{(1)}$ as an alternative realization of $\mathcal{N}^{(2)}$ or vice-versa.

In \cite{Sz2} the author presents an algorithm for producing sparse and dense alternative realizations of a given network $\mathcal{N}'$, i.e. realizations with the fewest and greatest number of reactions capable of generating the same kinetics (\ref{de}) as that given by $\mathcal{N}'$. Key to the analysis is fixing the matrix $Y$ to contain only the (source or product) complexes corresponding to the network $\mathcal{N}'$. The problem of finding an alternative realization $\mathcal{N}$ of $\mathcal{N}'$ then becomes one of finding a kinetics matrix $A_k$ such that
\[Y \cdot A_k \cdot \Psi(\mathbf{x}) = Y \cdot A_k' \cdot \Psi(\mathbf{x}).\]
If we set $M=Y \cdot A_k'$ and impose that $A_k$ be a kinetics matrix, dynamical equivalence can be guaranteed by the conditions
\begin{equation}
\label{realization}
\mbox{\textbf{(DE)}} \; \; \left\{ \; \; \begin{array}{ll} & \displaystyle{Y \cdot A_k = M} \\ & \displaystyle{\sum_{i=1}^m [A_k]_{ij} = 0, \; \; \; j=1, \ldots, m} \\ & \displaystyle{[A_k]_{ij} \geq 0, \; \; \; i,j = 1, \ldots, m, \; i \not= j} \\ & \displaystyle{[A_k]_{ii} \leq 0, \; \; \; i = 1, \ldots, m.} \end{array} \right.
\end{equation}

A sparse (respectively, dense) realization is given by a matrix $A_k$ satisfying (\ref{realization}) with the most (respectively, least) off-diagonal entries which are zeroes. A correspondence between the non-zero off-diagonal entries in $A_k$ and a positive integer value can be made by considering the binary variables $\delta_{ij} \in \left\{ 0, 1 \right\}$ which will keep track of whether a reaction is `on' or `off', i.e. we have
\[\delta_{ij} = 1 \leftrightarrow [A_k]_{ij} > \epsilon, \; \; \; i, j = 1, \ldots, m, \; \; i \not= j\]
for some sufficient small $0 < \epsilon \ll 1$, where the symbol `$\leftrightarrow$' denotes the logical relation `if and only if'. These proposition logic constraints for the structure of a network can then be formulated as the following linear mixed-integer constraints (see, for example, \cite{R-G}):
\begin{equation}
\label{density}
\mbox{\textbf{(S)}} \; \; \left\{ \; \; \begin{array}{ll} & \displaystyle{0 \leq [A_k]_{ij}-\epsilon \delta_{ij}, \; \; \; i,j = 1, \ldots, m, \; \; i \not= j} \\ & \displaystyle{0 \leq -[A_k]_{ij}+u_{ij} \delta_{ij}, \; \; \; i,j = 1, \ldots, m, \; \; i \not= j} \\ & \displaystyle{\delta_{ij} \in \left\{ 0, 1 \right\}, \; \; \; i,j = 1, \ldots, m, i \not= j,} \end{array} \right.
\end{equation}
where $u_{ij} > 0$ for $i, j = 1, \ldots, m, i \not= j$. The number of reactions present in the network corresponding to $A_k$ is then given by the sum of the $\delta_{ij}$'s so that the problem of determining a sparse network corresponds to satisfying the objective function
\begin{equation}
\label{sparse}
\mbox{\textbf{(Sparse)}} \; \; \left\{ \; \; \; \; \; \; \; \; \mbox{minimize} \; \; \; \; \; \sum_{i,j=1, i \not= j}^m \delta_{ij} \right.
\end{equation}
over the constraint sets (\ref{realization}) and (\ref{density}). Finding a dense network corresponds to maximizing the same function, which can also be stated as a minimization problem as
\begin{equation}
\label{dense}
\mbox{\textbf{(Dense)}} \; \; \left\{ \; \; \; \; \; \; \; \; \mbox{minimize} \; \; \; \; \; \sum_{i,j=1, i \not= j}^m -\delta_{ij}. \right.
\end{equation}

\subsection{Weakly Reversible Networks}
\label{weaklyreversiblenetworkssection}

Weakly reversible networks are a particularly important class of reaction networks because strong properties are known about their dynamics. Under a supplemental condition, which is easily derived from the reaction graph alone, it is known that there is a unique positive equilibrium concentration within each invariant space of the network and that that equilibrium concentration is at least locally asymptotically stable \cite{F1,H,H-J1}.

In \cite{Sz-H-T} the authors introduce an algorithm for determining dense weakly reversible realizations of a given kinetics. The algorithm is based on the fact that there are no cycles involving elements in different strongly connected components of a reaction network \cite{BJ-G}, and that for a fixed complex set the structure of the dense realization of a network is unique and contains the structures of all other possible realizations as sub-graphs\cite{Sz-H-P}. Omitting technical details, the algorithm can be summarized as:
\begin{enumerate}
\item
Set the matrices $Y$ and $M$ and initialize $\mathcal{K} = \left\{ \right\}$.
\item
Compute a dense realization $A_k$ forcing the edges in $\mathcal{K}$ to be excluded.
\item
Check whether $A_k$ is weakly reversible (if so, end algorithm and return $A_k$).
\item
Find all edges in $A_k$ which lead from one strongly connected component to another and add them to $\mathcal{K}$.
\item
Check whether these edges may be removed (if so, repeat steps (2)-(4); if not, end algorithm and return $A_k=0$).
\end{enumerate}

The algorithm has the drawbacks that it can only compute dense realizations and not sparse ones, and that it requires potentially multiple MILP optimizations which are known to be NP-hard. In Section \ref{weaklyreversible} we will present a method for determining both dense and sparse weakly reversible realizations in a single MILP optimization step.

\section{Original Results}

In this section, we extend the results of \cite{Sz-H-T} by showing how the requirement of weak reversibility can be formulated as a linear constraint. Consequently, the question of determining a sparse or dense weakly reversible realization can now be handled in a single MILP optimization step. We also extend this framework to cover the linearly conjugate networks introduced in \cite{J-S3}.

\subsection{Weak Reversibility as a Linear Constraint}
\label{weaklyreversible}

In this section we show that the requirement of weak reversibility can be formulated as a linear constraint. We require the following classical result about weakly reversible networks, which is modified from Theorem 3.1 of \cite{G-H} and Proposition 4.1 of \cite{F3}:
\begin{theorem}
\label{weaklyreversible}
Let $A_k$ be a kinetics matrix and let $\Lambda_i$, $i=1, \ldots, \ell,$ denote the support of the $i^{th}$ linkage class. Then the reaction graph corresponding to $A_k$ is weakly reversible if and only if there is a basis of ker$(A_k)$, $\left\{ \mathbf{b}^{(1)}, \ldots, \mathbf{b}^{(\ell)} \right\}$, such that, for $i=1, \ldots, \ell$,
\[\mathbf{b}^{(i)} = \left\{ \begin{array}{ll} b^{(i)}_j > 0, \hspace{0.3in} & j \in \Lambda_i \\ b^{(i)}_j = 0, & j \not\in \Lambda_i. \end{array} \right.\]
\end{theorem}

An immediate consequence of Theorem \ref{weaklyreversible} is that there is a vector $\mathbf{b} \in \mathbb{R}_{>0}^m \; \cap \; $ker$(A_k)$ if and only if the reaction graph corresponding to $A_k$ is weakly reversible. In other words, we can guarantee weak reversibility by imposing the condition
\begin{equation}
\label{231}
A_k \cdot \mathbf{b} = \mathbf{0}
\end{equation}
for some $\mathbf{b} \in \mathbb{R}_{>0}^m$. This is a nonlinear constraint in the $k_{ij}$'s and $b_j$'s. In order to make it linear, we consider the matrix $\tilde{A}_k$ with entries
\begin{equation}
\label{lala}
[\tilde{A}_k]_{ij} = [A_k]_{ij} \cdot b_j.
\end{equation}
It is clear from (\ref{lala}) that $\tilde{A}_k$ encodes a kinetics matrix and that $\mathbf{1} \in \mathbb{R}^m$ (the $m$-dimensional vector containing only ones) lies in ker$(\tilde{A}_k)$. Moreover, it is easy to see that $\tilde{A}_k$ encodes a weakly reversible network if and only if $A_k$ corresponds to a weakly reversible network. We can therefore check weak reversibility of the chemical reaction network corresponding to $A_k$ with the linear conditions
\begin{equation}
\label{weakreversibility1}
\mbox{\textbf{(WR')}} \; \; \left\{ \; \; \begin{array}{ll} & \displaystyle{\sum_{i=1}^m[\tilde{A}_k]_{ij} = 0, \; \; \; j=1, \ldots, m} \\ & \displaystyle{\sum_{i=1}^m[\tilde{A}_k]_{ji} = 0, \; \; \; j=1, \ldots, m} \\ & \displaystyle{[\tilde{A}_k]_{ij} \geq 0, \; \; \; i,j = 1, \ldots, m, \; i \not= j} \\ & \displaystyle{[\tilde{A}_k]_{ii} \leq 0, \; \; \; i = 1, \ldots, m.} \end{array} \right.
\end{equation}
By solving for the diagonal elements of $\tilde{A}_k$, the set of constraints (\ref{weakreversibility1}) can be simplified to
\begin{equation}
\label{weakreversibility}
\mbox{\textbf{(WR)}} \; \; \left\{ \; \; \begin{array}{ll} & \displaystyle{\sum_{i=1,i \not= j}^m[\tilde{A}_k]_{ij} = \sum_{i=1,i \not= j}^m[\tilde{A}_k]_{ji}, \; \; \; j=1, \ldots, m}\\ & \displaystyle{[\tilde{A}_k]_{ij} \geq 0, \; \; \; i,j = 1, \ldots, m, \; i \not= j.} \end{array} \right.
\end{equation}

No condition comparable to $Y \cdot A_k = M$ exists for the matrix $\tilde{A}_k$ so that we are left to optimization with respect to the internal entries of both $A_k$ and $\tilde{A}_k$. Given appropriate choices of $0 < \epsilon \ll 1$ and $u_{ij} > 0$, $i,j=1, \ldots, m$, $i \not= j$, we can impose
\begin{equation}
\label{density2}
\mbox{\textbf{(WR-S)}} \; \; \left\{ \; \; \begin{array}{ll} & \displaystyle{0 \leq [\tilde{A}_k]_{ij}-\epsilon \delta_{ij}, \; \; \; i,j = 1, \ldots, m, \; \; i \not= j}\\ & \displaystyle{0 \leq -[\tilde{A}_k]_{ij}+u_{ij} \delta_{ij}, \; \; \; i,j = 1, \ldots, m, \; \; i \not= j} \end{array} \right.
\end{equation}
as well as (\ref{density}) to ensure that both $A_k$ and $\tilde{A}_k$ contain zero and non-zero entries in the same places so that they correspond to reaction graphs with the same structure.

\subsection{Linear Conjugacy}
\label{linearconjugacysection}

In \cite{J-S2} the authors extended the concept of dynamical equivalence to linear conjugacy. In their framework, two networks $\mathcal{N}$ and $\mathcal{N}'$ are said to be linearly conjugate if there is a linear mapping which takes the flow of one network to the other. The case of two networks being realizations of the same kinetics is encompassed as a special case of linear conjugacy taking the transformation to be the identity. (For a more complete introduction to the notion of dynamical equivalence and conjugacy see \cite{P} or \cite{W}.)

Importantly, linearly conjugate networks share the same qualitative dynamics (e.g. number and stability of equilibria, persistence/extinction of species, dimensions of invariant spaces, etc.). Similarly with different realizations of the same kinetics (\ref{de}), if a network with unknown kinetics is linearly conjugate to a network with known dynamics, then the qualitative properties of the second network are transferred to the first.

The main result of \cite{J-S2} is the following. We have adopted the notation to match that contained in this paper. The notation is sufficiently distinct that we will prove the result independently here.
\begin{theorem}
\label{maintheorem}
Consider two mass-action systems $\mathcal{N} = (\mathcal{S},\mathcal{C},\mathcal{R})$ and $\mathcal{N}' = (\mathcal{S},\mathcal{C}',\mathcal{R}')$ and let $Y$ be the stoichiometric matrix corresponding to the complexes in either network. Consider a kinetics matrix $A_k$ corresponding to $\mathcal{N}$ and suppose that there is a kinetics matrix $A_b$ with the same structure as $\mathcal{N}'$ and a vector $\mathbf{c} \in \mathbb{R}_{>0}^n$ such that
\begin{equation}
\label{condition}
Y \cdot A_k = T \cdot Y \cdot A_b
\end{equation}
where $T = $diag$\left\{ \mathbf{c} \right\}$. Then $\mathcal{N}$ is linearly conjugate to $\mathcal{N}'$ with kinetics matrix
\begin{equation}
\label{newrateconstants}
A_k' = A_b \cdot \mbox{diag} \left\{ \Psi (\mathbf{c}) \right\}.
\end{equation}
\end{theorem}

\begin{proof}
Let $\Phi(\mathbf{x}_0,t)$ correspond to the flow of (\ref{de}) associated to the reaction network $\mathcal{N}$. Consider the linear mapping $\mathbf{h}(\mathbf{x})=T^{-1} \cdot \mathbf{x}$ where $T=$diag$\left\{ \mathbf{c} \right\}$. Now define $\tilde{\Phi}(\mathbf{y}_0,t)=T^{-1} \cdot \Phi(\mathbf{x}_0,t)$ so that $\Phi(\mathbf{x}_0,t) = T \cdot \tilde{\Phi}(\mathbf{y}_0,t)$.

Since $\Phi(\mathbf{x}_0,t)$ is a solution of (\ref{de}) we have
\[\begin{split} \tilde{\Phi}'(\mathbf{y}_0,t) & = T^{-1} \cdot \Phi'(\mathbf{x}_0,t) \\ & = T^{-1} \cdot Y \cdot A_k \cdot \Psi(\Phi(\mathbf{x}_0,t)) \\ & = T^{-1} \cdot T \cdot Y \cdot A_b \cdot \Psi(T \cdot \tilde{\Phi}(\mathbf{y}_0,t)) \\ & = Y \cdot A_b \cdot \mbox{diag} \left\{ \Psi (\mathbf{c}) \right\} \cdot \Psi (\tilde{\Phi}(\mathbf{y}_0,t)).\end{split}\]
\noindent It is clear that $\tilde{\Phi}(\mathbf{y}_0,t)$ is the flow of (\ref{de}) corresponding to the reaction network $\mathcal{N}'$ with the kinetics matrix given by (\ref{newrateconstants}). We have that $\mathbf{h}(\Phi(\mathbf{x}_0,t))=\tilde{\Phi}(\mathbf{h}(\mathbf{x}_0),t)$ for all $\mathbf{x}_0 \in \mathbb{R}_{>0}^n$ and $t \geq 0$ where $\mathbf{y}_0 = \mathbf{h}(\mathbf{x}_0)$ since $\mathbf{y}_0 = \tilde{\Phi}(\mathbf{y}_0,0)=T^{-1} \cdot \Phi(\mathbf{x}_0,0) = T^{-1} \cdot \mathbf{x}_0$. It follows that the networks $\mathcal{N}$ and $\mathcal{N}'$ are linearly conjugate and we are done.
\end{proof}

While the results of \cite{J-S2} give conditions for two networks to be linearly conjugate, and therefore exhibit the same qualitative dynamics, no general methodology is provided for determining linearly conjugate networks when only a single network is provided.

In cases where the dynamics of a network $\mathcal{N}$ is suspected to behave like a weakly reversible network, it is beneficial to extend the optimization algorithm outlined in Section \ref{weaklyreversible} to linearly conjugate networks. This can be accomplished by replacing the set of constraints (\ref{realization}) with
\begin{equation}
\label{conjugate}
\mbox{\textbf{(LC)}} \; \; \left\{ \; \; \begin{array}{ll} & \displaystyle{Y \cdot A_b = T^{-1} \cdot M} \\ & \displaystyle{\sum_{i=1}^m [A_b]_{ij} = 0, \; \; \; j=1, \ldots, m} \\ & \displaystyle{[A_b]_{ij} \geq 0, \; \; \; i,j = 1, \ldots, m, \; i \not= j} \\ & \displaystyle{[A_b]_{ii} \leq 0, \; \; \; i = 1, \ldots, m} \\ &\displaystyle{\epsilon  \leq c_j \leq 1/\epsilon, \; \; \; j=1, \ldots, n} \end{array} \right.
\end{equation}
where $M = Y \cdot A_k$, $T = $diag$\left\{ \mathbf{c} \right\}$, and $0 < \epsilon \ll 1$, and replacing the set of constraints (\ref{density}) by
\begin{equation}
\label{lcdensity}
\mbox{\textbf{(LC-S)}} \; \; \left\{ \; \; \begin{array}{ll} & \displaystyle{0 \leq [A_b]_{ij}-\epsilon \delta_{ij}, \; \; \; i,j = 1, \ldots, m, \; \; i \not= j} \\ & \displaystyle{0 \leq -[A_b]_{ij}+u_{ij} \delta_{ij}, \; \; \; i,j = 1, \ldots, m, \; \; i \not= j} \\ & \displaystyle{\delta_{ij} \in \left\{ 0, 1 \right\}, \; \; \; i,j = 1, \ldots, m, i \not= j,} \end{array} \right.
\end{equation}
where $u_{ij} > 0$ for $i, j = 1, \ldots, m, i \not= j$.

$A_b$ has the same structure as the kinetics matrix $A_k'$ corresponding to the conjugate network, and this matrix has the same structure as the matrix $\tilde{A}_k$ given by (\ref{lala}) (replacing $A_k$ by $A_k'$). Consequently, the problem of determining a sparse or dense weakly reversible network which is linearly conjugate to a given kinetics can be given by optimizing either (\ref{sparse}) or (\ref{dense}), respectively, over the constraint sets (\ref{conjugate}), (\ref{lcdensity}), (\ref{weakreversibility}), and (\ref{density2}). The kinetics matrix $A_k'$ for the linearly conjugate network is given by (\ref{newrateconstants}).

\section{Examples}

In this section we will consider two examples from the literature which demonstrate how the MILP optimization algorithm outlined in Section \ref{linearconjugacysection} is capable of efficiently finding sparse and dense weakly reversible networks which are linearly conjugate to a given network $\mathcal{N}$. We also consider one new example which illustrates how the algorithm is capable of finding networks with linearly conjugate dynamics for which no trivial linear conjugacy exists.\\

\textbf{Example 1:} Consider the chemical reaction network $\mathcal{N}$ given by
\[\begin{split} X_1 + 2X_2 \; & \stackrel{\alpha}{\longrightarrow} \; X_1 \\ \mathcal{N}: \hspace{0.5in} 2 X_1 + X_2 \; & \stackrel{1}{\longrightarrow} \; 3 X_2 \\ X_1 + 3 X_2 \; & \stackrel{1}{\longrightarrow} \; X_1 + X_2 \; \stackrel{1}{\longrightarrow} \; 3 X_1 + X_2. \hspace{0.3in} \end{split}\]

This network was first considered in \cite{J-S2} where the authors showed that it was linearly conjugate to a specified weakly reversible network for all values of $\alpha > 0$. It was further analysed with the value $\alpha = 1.5$ in \cite{Sz-H-T} where the authors found a dense weakly reversible realization through successive MILP optimizations. We will reproduce this result using our one-step MILP algorithm and also produce a sparse realization.

\begin{figure}[h]
\begin{center}
\includegraphics[width=11cm]{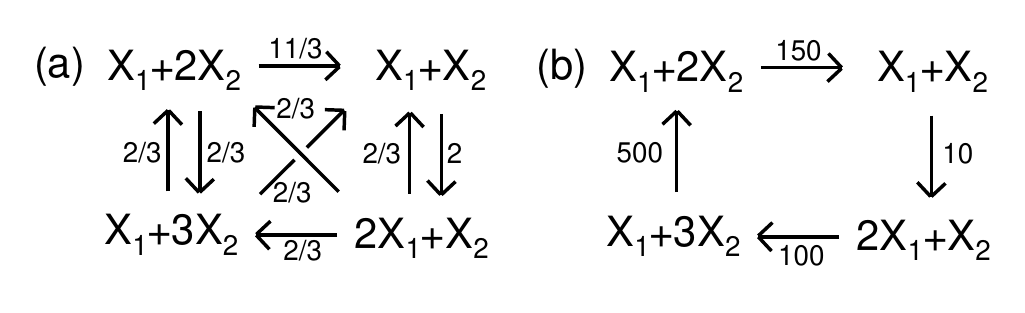}
\end{center}
\label{figure1}
\caption{Dense (a) and sparse (b) weakly reversible networks which are linearly conjugate to $\mathcal{N}$.}
\end{figure}

We have
\[Y = \left[ \begin{array}{ccccccc} 1 & 1 & 2 & 0 & 1 & 1 & 3 \\ 2 & 0 & 1 & 3 & 3 & 1 & 1 \end{array} \right] \]
and
\[M = \left[ \begin{array}{ccccccc} 0 & 0 & -2 & 0 & 0 & 2 & 0 \\ -3 & 0 & 2 & 0 & -2 & 0 & 0 \end{array} \right] \]
and set $\epsilon = 1 / \alpha = 2 / 3$ and $u_{ij} = 20$, $i,j = 1, \ldots, 7, i \not= j$. The MILP problem for a dense weakly reversible linearly conjugate network, possibly accounting for a non-trivial linear conjugacy mapping, is
\[\mbox{minimize} \; \; \; \; \; \sum_{i, j=1}^7 -\delta_{ij}\]
over the constraint set
\[\begin{split} Y \cdot A_b & = T^{-1} \cdot M \\ \sum_{i=1}^m[A_b]_{ij} & = 0, \; \; \; j=1, \ldots, m \\ \sum_{i=1,i \not= j}^m[\tilde{A}_k]_{ij} & = \sum_{i=1,i \not= j}^m[\tilde{A}_k]_{ji}, \; \; \; j=1, \ldots, m \\ 0 & \leq [A_b]_{ij}-\epsilon \cdot \delta_{ij}, i,j = 1, \ldots, 7, i \not= j \\ 0 & \leq -[A_b]_{ij} + u_{ij} \cdot \delta_{ij}, i,j = 1, \ldots, 7, i \not= j \\ 0 & \leq [\tilde{A}_k]_{ij}-\epsilon \cdot \delta_{ij}, i,j = 1, \ldots, 7, i \not= j \\ 0 & \leq -[\tilde{A}_k]_{ij} + u_{ij} \cdot \delta_{ij}, i,j = 1, \ldots, 7, i \not= j \end{split} \]
where $T = $diag$\left\{ \mathbf{c} \right\}$, and the decision variables
\[\begin{split} [A_b]_{ij} & \geq 0, [\tilde{A}_k]_{ij} \geq 0, \mbox{ for } i, j = 1, \ldots, 7, i \not= j \\ [A_b]_{ii} & \leq 0, \mbox{ for } i = 1, \ldots, 7 \\ \epsilon \leq c_i & \leq 1/\epsilon, \mbox{ for } i = 1, 2 \\ \delta_{ij} & \in \left\{ 0, 1 \right\}, \mbox{ for } i,j = 1, \ldots, 7, i \not= j.\end{split}\]
Solving for $A_b$ with GLPK and applying (\ref{newrateconstants}) gives the kinetics matrix
\[A_k = \left[ \begin{array}{ccccccc} -\frac{13}{3} & 0 & \frac{2}{3} & 0 & \frac{2}{3} & 0 & 0 \\ 0 & 0 & 0 & 0 & 0 & 0 & 0 \\ 0 & 0 & -2 & 0 & 0 & 2 & 0 \\ 0 & 0 & 0 & 0 & 0 & 0 & 0 \\ \frac{2}{3} & 0 & \frac{2}{3} & 0 & -\frac{4}{3} & 0 & 0 \\ \frac{11}{3} & 0 & \frac{2}{3} & 0 & \frac{2}{3} & -2 & 0 \\ 0 & 0 & 0 & 0 & 0 & 0 & 0 \end{array} \right]\]
and values $c_1 = 1, c_2 = 1$ (i.e. the linear transformation is the identity). The network structure is given graphically in Figure \ref{figure1}(a). Although the rate constants differ due to differing bounds, this has the same network structure as the dense weakly reversible network obtained in \cite{Sz-H-T}.

A sparse weakly reversible network is generated by optimizing
\[\mbox{minimize} \; \; \; \; \; \sum_{i, j=1}^7 \delta_{ij}\]
over the same constraint set.  Solving for $A_b$ with GLPK with the bound $\epsilon=0.1$ and applying (\ref{newrateconstants}) gives the kinetic matrix
\[A_k = \left[ \begin{array}{ccccccc} -150 & 0 & 0 & 0 & 500 & 0 & 0 \\ 0 & 0 & 0 & 0 & 0 & 0 & 0 \\ 0 & 0 & -100 & 0 & 0 & 10 & 0 \\ 0 & 0 & 0 & 0 & 0 & 0 & 0 \\ 0 & 0 & 100 & 0 & -500 & 0 & 0 \\ 150 & 0 & 0 & 0 & 0 & -10 & 0 \\ 0 & 0 & 0 & 0 & 0 & 0 & 0 \end{array} \right]\]
and values $c_1 = 10$ and $c_2 = 5$. This is therefore an example of a network with a non-trivial linear conjugacy and corresponds to the weakly reversible network given in Figure \ref{figure1}(b).\\

\textbf{Example 2:} Consider the kinetics scheme
\begin{equation}
\label{example}
\begin{split} \dot{x}_1 & = x_3^2 - x_1 x_2 + x_3 x_4 - 2 x_1 x_2^2 x_3 \\ \dot{x}_2 & = x_3^2 - x_1 x_2 + 2 x_3 x_4 - 4 x_1 x_2^2 x_3 \\ \dot{x}_3 & = -2 x_3^2 + x_1 x_2 - x_1 x_2^2 x_3 + 2 x_4^3 \\ \dot{x}_4 & = x_1 x_2 - x_3 x_4 + 4 x_1 x_2^2 x_3 - 3 x_4^3 \end{split}
\end{equation}
first considered in \cite{Sz-H}. Using the algorithm given in \cite{H-T} and \cite{Sz-H}, we can determine a kinetic realization involving the complexes
\[\begin{split} & C_1 = 2 X_3, C_2 = X_3 + X_4, C_3 = X_1 + 2 X_3, C_4 = X_2 + 2 X_3, \\ & C_5 = X_3, C_6 = X_1 + X_3 + X_4, C_7 = X_2 + X_3 + X_4 \\ & C_8 = X_1 + X_2, C_9 = X_1 + 2 X_2 + X_3, C_{10} = X_1, C_{11} = X_2 \\ & C_{12} = X_1 + X_2 + X_4, C_{13} = X_1 + X_2 + X_3, C_{14} = 2 X_2 + X_3 \\ & C_{15} = X_1 + 2 X_2, C_{16} = X_1 + 2 X_2 + X_3 + X_4\\ & C_{17} = 3 X_4, C_{18} = X_3 + 3 X_4, C_{19} = 2 X_4.\end{split}\]

\begin{figure}[h]
\begin{center}
\includegraphics[width=8cm]{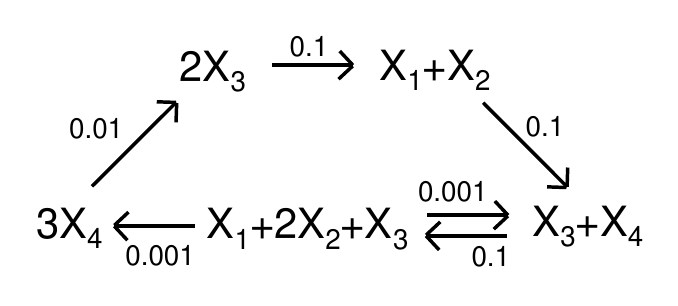}
\end{center}
\label{figure2}
\caption{Weakly reversible realization of the kinetics (\ref{example}). This realization is both dense and sparse.}
\end{figure}

In \cite{Sz-H-T} the authors use the algorithm given in Section \ref{weaklyreversiblenetworkssection} to determine a dense weakly reversible realization. The algorithm required three MILP optimizations, three searches for strongly connected components, and took 80.5s to complete. Carrying out either the MILP optimization algorithm outlined in Section \ref{linearconjugacysection} for a dense or for a sparse weakly reversible network, with bounds $\epsilon = 0.1$ and $u_{ij} = 10,$ $i,j=1, \ldots, 19$, we arrive at the solution
\[\tilde{k}_{18} = \tilde{k}_{29} = \tilde{k}_{82} = 0.1, \; \; \tilde{k}_{92} = \tilde{k}_{9(17)} = 0.001, \; \; \tilde{k}_{(17)1} = 0.01,\]
\[c_1 =c_2=c_3=c_4= 0.1\]
and the rest of the entries zero (the transformation is a scaling of the identity). This corresponds to the network given in Figure \ref{figure2} which has the same network structure as the networks obtained in \cite{Sz-H} and \cite{Sz-H-T}. Our algorithm was able to obtain the answer in a single MILP optimization step and took less than a tenth of a second to compute. (The difference in rate constants occurs as a result of the scaling of concentration variables permitted by linear conjugacy.)\\

\textbf{Example 3:} Consider the kinetics scheme
\begin{equation}
\label{example2}
\begin{split} \dot{x}_1 & = x_1 x_2^2 - 2 x_1^2 + x_1 x_3^2 \\ \dot{x}_2 & = -x_1^2 x_2^2 + x_1 x_3^2 \\ \dot{x}_3 & = x_1^2 - 3 x_1 x_3^2. \end{split}
\end{equation}
Using the algorithm given in \cite{H-T} and \cite{Sz-H}, we can determine a kinetic realization involving the complexes
\[\begin{split} & C_1 = X_1 + 2 X_2, C_2 = 2X_1 + 2 X_2, C_3 = 2X_1 + X_2, \\ & C_4 = 2X_1, C_5 = X_1, C_6 = 2X_1 + X_3, C_7 = X_1 + 2X_3 \\ & C_8 = 2X_1 + 2X_3, C_9 = X_1 + X_2 + 2X_3, C_{10} = X_1 + X_3. \end{split}\]

With this fixed complex set, we can carry out the MILP optimization procedure outlined in Section \ref{linearconjugacysection} to find sparse and dense weakly reversible networks which are linearly conjugate to a network with kinetics (\ref{example2}). We have
\[Y= \left[ \begin{array}{cccccccccc} 1 & 2 & 2 & 2 & 1 & 2 & 1 & 2 & 1 & 1 \\ 2 & 2 & 1 & 0 & 0 & 0 & 0 & 0 & 1 & 0 \\ 0 & 0 & 0 & 0 & 0 & 1 & 2 & 2 & 2 & 1 \end{array} \right] \]
and
\[M = \left[ \begin{array}{cccccccccc} 1 & 0 & 0 & -2 & 0 & 0 & 1 & 0 & 0 & 0 \\ 0 & -1 & 0 & 0 & 0 & 0 & 1 & 0 & 0 & 0 \\ 0 & 0 & 0 & 1 & 0 & 0 & -3 & 0 & 0 & 0 \end{array} \right].\]

\begin{figure}[h]
\begin{center}
\includegraphics[width=13cm]{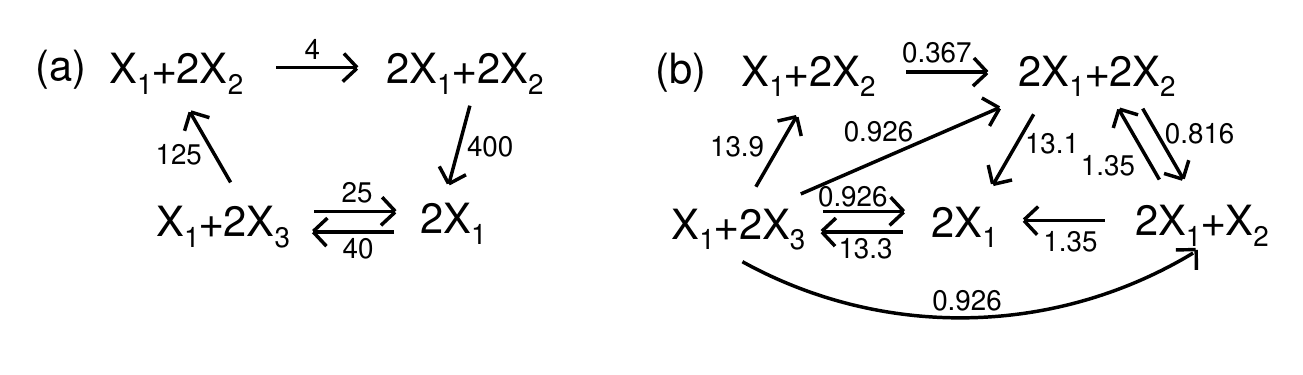}
\end{center}
\label{figure3}
\caption{Weakly reversible networks which are linearly conjugate to a network with the kinetics (\ref{example2}). The network in (a) is sparse while the network in (b) is dense. The parameter values in (b) have been rounded to three significant figures.}
\end{figure}

With the bounds $\epsilon=1/20$ and $u_{ij}=20$ for $i,j=1, \ldots, 10,$ $i \not= j$, the algorithm gives us the sparse network given in Figure \ref{figure3}(a) (conjugacy constants $c_1=20$, $c_2=2$, and $c_3=5$) and the dense network given in Figure \ref{figure3}(b) (conjugacy constants $c_1=20/3$, $c_2=20/33$, and $c_3=5/3$). It is interesting to note that the sparse and dense networks utilize different complexes and that the ratio of conjugacy constants differ between the sparse and dense networks. It is worth noting that the sparse realization is also deficiency zero so that the Deficiency Zero Theorem can be applied \cite{F1,H,H-J1}. Consequently, solutions of (\ref{example2}) satisfy all of the stringent dynamical restrictions typically reserved for complex balanced systems.

\section{Conclusions}

In this paper we have proposed an algorithm for determining linearly conjugate weakly reversible realizations of reaction networks. In contrast to the method presented in \cite{Sz-H-T}, the present approach is based on the well-known fact that the kernel of the Kirchhoff matrix of weakly reversible networks always contains a strictly positive vector. The main advantages of our algorithm compared to \cite{Sz-H-T} are the following. Firstly, linear conjugacy theory \cite{J-S2} has been included into the optimization framework, and the parameters of the corresponding linear coordinates transformation belong to the set of unknowns. Secondly, our algorithm requires only one MILP step, and therefore it is numerically significantly more effective than \cite{Sz-H-T} if the problem dimension is similar to what is shown in the examples. Thirdly, additional structural constraints such as density or sparsity of the solution network can be directly included into the optimization problem. The presented results clearly contribute to further widening the application possibilities of the known strong results in chemical reaction network theory.

There are still several very important questions which remain to be answered:
\begin{enumerate}
\item
While the algorithm is effective and efficient for finding alternative networks for a given kinetics, we are often interested in questions which can be answered for all kinetics satisfying certain initial structural properties (i.e. general rather than specied rate constants). The algorithm is currently unable to answer such questions.
\item
Many dynamical properties are known for systems satisfying network structure properties not included in weak reversibility theory \cite{C-F1,C-F2,C-P2}. Extending the optimization framework to include these results would greatly expand the scope of chemical reaction networks with known dynamics.
\end{enumerate}

\section*{Acknowledgements}

M. Johnston and D. Siegel acknowledge the support of D. Siegel's Natural Sciences and Engineering Research Council of Canada Discovery Grant.\\

\noindent G. Szederk\'{e}nyi acknowledges the support of the Hungarian National Research Fund through grant no. OTKA K-83440 as well as the support of project CAFE (Computer Aided Process for Food Engineering)  FP7-KBBE-2007-1
(Grant no: 212754).

\bibliographystyle{plain}
\bibliography{myrefs}

.

\end{document}